\documentclass[12pt]{article}

\usepackage{amsfonts}
\usepackage{amsthm}
\usepackage[mathscr]{eucal}
\usepackage{amsmath}
\usepackage{amssymb}
\usepackage{amscd}

\theoremstyle{plain}
\newtheorem{lemma}{Lemma}[section]
\newtheorem{theorem}[lemma]{Theorem}
\newtheorem{proposition}[lemma]{Proposition}
\newtheorem{corollary}[lemma]{Corollary}

\theoremstyle{definition}

\newtheorem{definition}[lemma]{Definition}
\newtheorem{remark}[lemma]{Remark}

\numberwithin{equation}{section}
\def\begeq{\stepcounter{lemma}\begin{equation}}

\date{}

\begin{document}

\title{The Other Group of a Galois Extension\\
(Castor and Pollux)\footnote{Castor and Pollux are the two stars that give the
constellation Gemini (Latin for ``twins") its name. In reality these two celestial
entities are unlikely twins. Castor is a hot, white system of six stars, while Pollux
is a single, large, cooler, yellow-orange star. The Chinese recognized Castor and Pollux
as Yin and Yang, the two fundamental, complementary forces of the
universe.} }

\author{ Lex E. Renner}
\maketitle

\begin{quote}
    {\em ``Both buried now in the fertile earth though still alive.
     Even under the earth Zeus grants them this distinction:
     one day alive, the next day dead, each twin takes his turn
     to have status among the gods."}  \;\;\;  The Odyssey\\
\end{quote}


\begin{abstract}

\noindent Let $k\subseteq K$ be a finite Galois extension of fields with Galois group $G$.
Let $\mathscr{G}$ be the automorphism $k$-group scheme of $K$.
We construct a canonical $k$-subgroup scheme $\underline{G}\subset\mathscr{G}$
with the property that $Spec_k(K)$ is a $k$-torsor for $\underline{G}$.
$\underline{G}$ is a constant $k$-group if and only if $G$ is abelian, in which case
$G=\underline{G}$.
\end{abstract}


\section{Introduction}   \label{intro.sec}

\vspace{10pt}

Let $k\subset K$ be a finite Galois extension of degree $n>2$
and let $\Gamma=Gal(k_s/k)$. Consider
the automorphism group $\mathscr{G}=Aut(K/k)$ as a functor from $k$-algebras
to $\Gamma$-groups,
\[
L\leadsto\mathscr{G}(L)=Aut_L(L\otimes_kK).
\]
As such $\mathscr{G}$ is representable by a finite \'etale $k$-group. It
turns out that $\mathscr{G}$ can be described by a homomorphism
\[
\theta : \Gamma\to Aut(S_n)
\]
such that $\theta(\Gamma)\subset Inn(S_n)$ is a simply transitive subgroup.
We can therefore write $\theta(\sigma)(g)=\beta(\sigma)g\beta(\sigma)^{-1}$ where
$\beta : \Gamma\to S_n$ is a group homomorphism.
It turns out that
\[
G=S_n^\Gamma=\{g\in S_n\;|\;\beta(\sigma)g=g\beta(\sigma)\;\text{for
all}\;\sigma\in\Gamma\}
\]
and
\[
\beta(\Gamma)=C_{S_n}(G)=\{h\in S_n\;|\;hg=gh\;\text{for all}\; g\in G\}
\]
where $G=Aut_k(K)=\mathscr{G}(k)\subset\mathscr{G}(K)=S_n$.
This yields an \'etale $k$-group $(\underline{G},\theta)$ by restricting $\theta$
to obtain $\theta|C_{S_n}(G) : \Gamma\to Aut(C_{S_n}(G))$, via
\[
\theta(\sigma)(g)=\beta(\sigma)g\beta(\sigma)^{-1}.
\]
for $g\in C_{S_n}(G)$ and $\sigma\in\Gamma$.
$(\underline{G},\theta)$ is thus a $k$-subgroup of $(S_n,\theta)$.

\begin{theorem}  \label{mainthm.thm}
Let $X=Spec_k(K)$. The restriction
\[
\underline{G}\times_k X\to X
\]
makes $X$ into a principal homogeneous space for $\underline{G}$.
\end{theorem}

\begin{theorem}   \label{maincor.cor}
$\underline{G}$ is a constant $k$-group if and only if
$G$ is abelian.
\end{theorem}


\section{Simply Transitive Groups}     \label{one.sec}

In this section we assemble a few basic results about finite groups
that are needed for our main results.

Consider the finite group $G$ of order $n$ thought of as a group of transformations
acting on itself by either left translations or right translations.
We write $l_g(x)=gx$ and $r_g(x)=xg$. If we think of $G$ acting on
itself by left translations we write $G^l=\{l_g\}\subseteq S_n$.
If we think of $G$ acting on
itself by right translations we write $G^r=\{r_g\}\subseteq S_n$.
We write $G$ if we think of $G$ as the set being acted on somehow,
or if we think of $G$ acting somehow other than on itself by right or
left translation. Notice that $IG^lI^{-1}=G^r$ where $I\in S_n$ is defined
by $I(g)=g^{-1}$. Indeed, for any $g\in G$, $I(l_g)I^{-1}=r_{g^{-1}}$.

\begin{proposition}   \label{cent.prop}
$C_{S_n}(G^l)=G^r$ and $C_{S_n}(G^r)=G^l$. Furthermore $G^l\cap G^r\cong Z(G)$,
the center of $G$.
\end{proposition}
\begin{proof}
It is straightforward to check that $G^r\subseteq C_{S_n}(G^l)$. Suppose conversely that
$f\in C_{S_n}(G^l)$ and let $g=f(1)\in G$. Then, for $x\in G$,
$f(x)=f(x1)=(f\circ l_x)(1)=(l_x\circ f)(1)=xf(1)=xg$. Evidently
$f=r_g\in G^r$.

Suppose that $f\in G^l\cap G^r$. Then $f=l_g=r_h$ for $g,h\in G$. But then
$f(1)=g=l_g(1)=r_h(1)=h$. On the other hand, for any $x\in G$,
$gx=l_g(x)=f(x)=r_h(x)=xh$. Conversely, it is trivial to check that
$Z(G^r)\cup Z(G^l)\subset G^r\cap G^l$.
\end{proof}

The following Proposition will be useful in the next section.

\begin{proposition}  \label{bigcentral.prop}
Suppose that $H<G^r$ is a subgroup of order $k$ so that $n=km$.
Then $|C_{S_n}(H)|=(m!)k^m$.
\end{proposition}
\begin{proof}
$H=\{r_g\;|\; g\in H\}$. So let
\[
G_r=\bigsqcup_{i=1}^mx_iH=\bigsqcup_{i=1}^mA_i.
\]
Let
\[
\sigma\in S_n\;\text{and}\; y_i\in A_i, i=1,...,m.\;\;\;(*)
\]
Define a bijection $f_* : G\to G$ by setting $f_*(x_ih)=y_{\sigma(i)}h$
so that $f_*(A_i)=A_{\sigma(i)}$. Then, for $g\in H$,
\[
f_*\circ r_g(x_ih)=f_*(x_ihg)=y_{\sigma(i)}hg=r_g(y_{\sigma(i)}h)=r_g\circ
f_*(x_ih).
\]

Conversely, given $f\in C_{S_n}(H)$, let $f(x_i)=y_{\sigma(i)}\in
A_{\sigma(i)}$. One checks that if $i\neq j$ then $\sigma(i)\neq
\sigma(j)$, since $f$ is injective. But then $f(x_ih)=y_{\sigma(i)}h$
where $\sigma(i)\in S_m$. Thus $f$ is among those already listed in $(*)$.
Counting these all up we obtain that there are $(m!)k^m$ different possibilities.
\end{proof}

\begin{corollary}        \label{bigcentral.cor}
Suppose that $m>1$ and $n>2$. Let $H<G^r$ be such that $C_{S_n}(H)=G^l$. Then $H=G^r$.
\end{corollary}
\begin{proof}
If $m>1$ and $n>2$ then $(m!)k^m>n$.
\end{proof}

\begin{definition}    \label{simplytrans.def}
A subgroup $G\subset S_n$ is called {\em simply transitive} if $G$ acts
simply and transitively on the set $\{1,2,...,n\}$ for some transitive action
of $S_n$ on $\{1,2,...,n\}$.
\end{definition}

\begin{remark}  \label{simplytranssaywat.rk}
One might wonder in what sense this transitive action is well-defined. If $n=2$ then $G=S_n$.
So that case is alright. If $n\neq 6$ then the action of $S_n$ on $\{1,2,...,n\}$ is well-defined
up to inner automorphism, so this case is also fine; one just takes the action of $G$ to be
conjugation on the $n$ conjugate subgroups of $S_n$ of index $n$. However there is a nontrivial
outer automorphism of $S_6$. See \cite{Lor}. Correspondingly, there are two conjugacy classes
of subgroups of $S_6$ isomorphic to $S_5$ and, consequently, there are two different ways to define
the cycle structure on the elements of $S_6$. Let $G\subset S_6$ be a subgroup that is simply
transitive on $\{1,2,3,4,5,6\}$ for the transitive action $T_1$ of $S_6$ on $\{1,2,3,4,5,6\}$.
Now $G$ is isomorphic to $\mathbb{Z}/6\mathbb{Z}$ or $S_3$. Let $g\in G$ be an element of order
two. Since $G$ is simply transitive for the $T_1$ action, $g$ has cycle structure $(2,2,2)$
for this action. But it is known (see \cite{Lor}) that, for the ``other" transitive action
$T_2$ of $S_6$ on $\{1,2,3,4,5,6\}$, $g$ must have cycle structure $(1,1,1,1,2)$. So $g$ has
four fixed points for the $T_2$ action. We conclude that, if $G\subset S_6$ is simply transitive
for one structure, then it is not simply transitive for the other (nonconjugate) structure.
\end{remark}


\section{Finite \'Etale $k$-groups}       \label{two.sec}

Let $k$ be  a field and let $\Gamma$ be the Galois group of the separable closure
$k_s$ over $k$. We first remind the reader of the well-known equivalence
between finite $\Gamma$-sets and finite separable $k$-algebras.

If $\Gamma\times X\to X$ is a $\Gamma$-set we denote the action
of $\sigma\in\Gamma$ on $x\in X$ by $\sigma(x)$. If $\Gamma\times X\to X$
is a finite $\Gamma$-set we define
\[
A_X=Hom_s(X,k_s)^\Gamma
\]
where the action of $\Gamma$ on $Hom_s(X,k_s)$ is given by
$\sigma(f)(x)=\sigma(f(\sigma^{-1}(x)))$. Here $\sigma$ also acts on $k_s$
and ``$Hom_s(U,V)$" denotes the set of functions from $U$ to $V$.
One checks that $A_X$ is a finite separable $k$-algebra.

If $k\to A$ is a finite separable $k$-algebra we define
\[
X_A=Hom_k(A,k_s)
\]
where the action of $\Gamma$ on $X_A$ is given by its action on $k_s$.
i.e. $\sigma(x)(a)=\sigma(x(a))$. ``$Hom_k(A,B)$" denotes the set of
$k$-algebra homorphisms from $A$ to $B$.

There is a canonical isomorphism of $\Gamma$-algebras,
\[
A\otimes_k k_s\cong Hom_s(X_A,k_s),
\]
defined by $a\otimes\lambda\leadsto f$, where $f(x)=x(a)\lambda$.

\begin{theorem}   \label{ksepvsgammasets.thm}
The functors $A\leadsto A_X$ and $X\leadsto X_A$ determine an
equivalence between the category of finite separable $k$-algebras
and the category of finite $\Gamma$-sets.
\end{theorem}
\begin{proof}
See Theorem 6.4 of \cite{Water}.
\end{proof}

With this in mind we make the following definition.

\begin{definition}         \label{finiteetalegroup.def}
A {\em finite \'etale $k$-group} $(\underline{G},\beta)$ is a finite abstract group $G$ together
with a homomorphism of groups $\beta : \Gamma\to Aut(G)$. The $k$-Hopf algebra of functions on
$\underline{G}$ is $k[\underline{G}]=\{f:G\to k_s\;|\; f(\beta(\sigma)(g))=\beta(\sigma)(f(g))\;\text{for
all}\;\sigma\in\Gamma\}$. $G$ is called a {\em $\Gamma$-group}.
\end{definition}

Using Theorem~\ref{ksepvsgammasets.thm} we obtain the following result.

\begin{theorem}   \label{twistedkgroups.thm}
There is a canonical equivalence between the category of finite
$\Gamma$-groups and the category of finite separable $k$-Hopf algebras.
\end{theorem}

\section{The Other Group}     \label{three.sec}

Let $k\subset K$ be a finite Galois extension of degree $n$ with Galois group $G$,
and let $\Gamma = Gal(k_s/k)$ be the Galois group of the separable closure $k_s$
of $k$. Define a finite \'etale $k$-group by the rule
\[
\mathscr{G}(L)=Aut_L(L\otimes_k K)
\]
where $k\subset L$ is a $k$-algebra. Notice that $\mathscr{G}(k_s)\cong S_n$.
By Theorem~\ref{twistedkgroups.thm} there is a homorphism of groups
\[
\theta : \Gamma \to Aut(S_n)
\]
such that $S_n^\Gamma=G\subset S_n$. Furthermore, $G\subset S_n$ is a simply transitive
subgroup in the sense of Definition~\ref{simplytrans.def}.

\begin{proposition}    \label{thetaisinner.prop}
Assume $n>2$. Then the group homomorphism
\[
\theta : \Gamma\to Aut(S_n)
\]
factors through $int : S_n\to Aut(S_n)$. Thus $\Gamma$ acts on $S_n$ by
\[
g\leadsto\beta(\sigma)g\beta(\sigma)^{-1}
\]
where $\beta : \Gamma\to S_n$ is a group homomorphism.
\end{proposition}
\begin{proof}
There is no doubt here unless $n=6$. Now $G\subset S_6$, being simply transitive, must be
isomorphic to $\mathbb{Z}/6\mathbb{Z}$ or $S_3$, and in such a way that any element
$g\in G=S_6^\Gamma$ of order two is conjugate in $S_6$ to $(12)(34)(56)$ (see
Remark~\ref{simplytrans.def}). But it is known \cite{Lor} that any outer automorphism of $S_6$
exchanges the conjugacy class of $(12)(34)(56)$ with that of $(12)$. Thus, if $\sigma\in\Gamma$
and $\theta(\sigma)$ is outer then $\theta(\sigma)(g)\neq g$. Thus $\theta(\sigma)$ is inner
for all $\sigma\in\Gamma$.
\end{proof}

\begin{corollary}    \label{itsinner.cor}
Assume $n>2$. Suppose that $\theta : \Gamma\to Aut(S_n)$ is such that $G=S_n^\Gamma$ is a
simply transitive subgroup. Then
\begin{enumerate}
\item $\theta(\Gamma)\subseteq Inn(S_n)$.
\item $\theta(\Gamma)=C_{S_n}(G)$.
\end{enumerate}
\end{corollary}
\begin{proof}
Corollary~\ref{bigcentral.cor} and Proposition~\ref{thetaisinner.prop}.
\end{proof}

Let $k\subset K$ be a finite Galois extension of degree $n>2$
and let $\Gamma=Gal(k_s/k)$. Consider
the Galois group $\mathscr{G}=Gal(K/k)$ as a functor from $k$-algebras
to $\Gamma$-groups ,
\[
L\leadsto\mathscr{G}(L)=Aut_L(L\otimes_kK).
\]
As such $\mathscr{G}$ is representable by a finite \'etale $k$-group.
As above $\mathscr{G}$ can be described by a homomorphism
\[
\theta : \Gamma\to Aut(S_n)
\]
such that $\theta(\Gamma)\subset Inn(S_n)$. We can therefore write
$\theta(\sigma)(g)=\beta(\sigma)g\beta(\sigma)^{-1}$ where
$\beta : \Gamma\to S_n$ is a group homomorphism. Furthermore
$\beta(\Gamma)\subset S_n$ is a simply transitive subgroup.
Now,
\[
G=S_n^\Gamma=\{g\in S_n\;|\;\beta(\sigma)g=g\beta(\sigma)\;\text{for
all}\;\sigma\in\Gamma\}
\]
and
\[
\theta(\Gamma)=C_{S_n}(G)=\{h\in S_n\;|\;hg=gh\;\text{for all}\; g\in G\}.
\]
where $G=Aut_k(K)=\mathscr{G}(k)\subset\mathscr{G}(K)=S_n$.
This yields an \'etale $k$-group $(\underline{G},\theta)$ by restricting $\theta$
to obtain $\theta|C_{S_n}(G) : \Gamma\to Aut(C_{S_n}(G))$, via
\[
\theta(\sigma)(g)=\beta(\sigma)g\beta(\sigma)^{-1}.
\]
for $g\in C_{S_n}(G)$.
$(\underline{G},\theta)$ is thus a $k$-subgroup of $(S_n,\theta)$. It is the centralizer
in $(S_n,\theta)$ of $(G,\theta)$, the latter being a constant subgroup since
$G=S_n^\Gamma$.

\begin{theorem}  \label{mainthm2.thm}
Let $X=Spec_k(K)$. The restriction
\[
\underline{G}\times_k X\to X
\]
makes $X$ into a principal homogeneous space for $\underline{G}$.
\end{theorem}
\begin{proof}
$\underline{G}(K)\times X(K)\to X(K)$ acts simply and
transitively by Proposition~\ref{cent.prop}. Thus the map
$\underline{G}(K)\times X(K)\to X(K)\times X(K)$, $(g,x)\to (gx,x)$,
is a bijection.
\end{proof}

\begin{theorem}   \label{maincor2.cor}
$\underline{G}\cap G = Z(G)$, the center of $G$.
$\underline{G}$ is a constant $k$-group if and only if
$G$ is abelian. In this case $\underline{G}=G$.
\end{theorem}
\begin{proof}
See Proposition~\ref{cent.prop}.
\end{proof}




\vspace{20pt}

\noindent Lex E. Renner \\
Department of Mathematics \\
University of Western Ontario \\
London, N6A 5B7, Canada \\

\enddocument